\documentclass[12pt,a4paper]{amsart}
\usepackage[utf8]{inputenc}
\usepackage{amsmath}
\usepackage{amsfonts}
\usepackage{amssymb}

\usepackage{hyperref}

\usepackage{float}
\usepackage{subfig}

\usepackage{graphicx}

\usepackage{caption}

\newtheorem{thm}{Theorem}[section]

\def\HH{\mathbb{H}}

\def\ZZ{\mathbb{Z}}

\def\RR{\mathbb{R}}
\def\QQ{\mathbb{Q}}
\def\NN{\mathbb{N}}

\def\tt{\Sigma_{1,1}}

\def\F2{\ZZ * \ZZ}
\def\aut{\text{Aut}(\F2)}
\def\gl2{\mathrm{GL}(2, \ZZ)}
\def\sl2{\mathrm{SL}(2, \RR)}

\def\hc{\tilde{\Sigma}}
\def\isom{\mathrm{isom}(\HH)}

\def\tr{\text{tr\,}}

\title{Convexity and Aigner's Conjectures}

 \author[McShane]{Greg McShane}
\address{Institut Fourier 100 rue des maths, BP 74, 38402 St Martin d'H\`eres cedex, France}
\email{mcshane at ujf-grenoble.fr}

\begin{document}

\maketitle

\begin{abstract} 

Markov numbers are integers that appear in triples which are solutions of a Diophantine equation, 
the so-called Markov cubic
 $$x^2 + y^2 + z^2 - 3x y z = 0.$$
 A classical topic in number
theory, these numbers are related to many areas of mathematics such as combinatorics,
hyperbolic geometry, approximation theory and cluster algebras.

One can associate  to each
a positive rational number a Markov number
in a natural way. 
We give a new  unified proof of certain conjectures from
Martin Aigner’s book, Markov’s
Theorem and 100 Years of the Uniqueness Conjecture. 
Our proof relies on a relationship
between Markov numbers and the lengths
of closed simple geodesics on the 
punctured torus discovered by H. Cohn.
\end{abstract}

\section{Introduction}

The Markov numbers are the positive integer solutions of the Diophantine equation 

$$ x^2 + y^2 + z^2 - 3xyz = 0 $$

Markov in the late 1800s
showed that all these solutions 
could be given the structure of a binary tree. 
It became customary (and useful) to index the Markov numbers by the rationals 
in the interval  $[0, 1]$
which stand at the same place in the
Stern–Brocot binary tree
(see Figures \ref{markov tree} and \ref{farey tree}). 
Frobenius’ conjecture asserts that each
Markov number appears at most once in this tree.
If this conjecture is true,
then the obvious order of Markov numbers
gives rise to a new strict order on the rationals.

\subsection{Aigner's Monotonicity Conjectures}
Aigner \cite{ai} proposed 
three conjectures 
to better understand this order.

\begin{enumerate}

\item (The fixed numerator conjecture) 
Let $p, q$ and $i$ be positive integers
such that $i > 0, p < q$
and $gcd(p, q) = 1,  gcd(p, q + i) = 1$ then 
$$m_{p/q} < m_{p/(q+i)}.$$

\item (The fixed denominator conjecture) 
Let $p, q$ and $i$ be positive integers
such that $i > 0, p < q$
and $gcd(p, q) = 1,  gcd(p + i,q) = 1$ then 
$$m_{p/q} < m_{(p+i)/q}.$$

\item (The fixed sum conjecture) 
Let $p, q$ and $i$ be positive integers
such that $i > 0, p < q$
and $gcd(p, q) = 1,  gcd(p - i, q + i) = 1$ then 
$$m_{p/q} < m_{(p-i) /(q+i)}.$$

\end{enumerate}

\begin{figure}[ht]
\centering
\includegraphics[scale=.7]{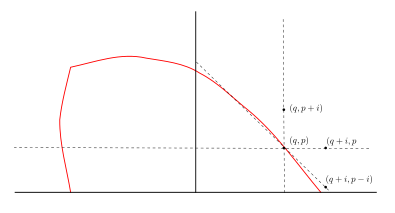} 
  \caption{A portion of the set $B$ with the point $(q,p)$ 
  on its boundary $S$ and the three lines 
  (horizontal, vertical and diagonal) 
  that appear in our proof of the conjectures.}
  \label{proof diagram}
\end{figure}

The fixed numerator conjecture was proved in \cite{ra} using ideas from cluster algebra theory
and snake graphs.
Their proof relies on a connection between Christoffel paths and Markov numbers (see \cite{berst} for a discussion)
which provides a formulation of Markov numbers as continuant polynomials.
The proof is rather technical 
and introduces several interesting tools 
for studying paths 
in the Cayley graph 
of the monoid  $\NN^2$ 
with the usual generators. 
In \cite{vu}, the authors prove the fixed numerator conjecture by using
"transformations" on paths in the Cayley graph.
We will show how these results can
be deduced by considerations of convexity
established in \cite{mcr}.
Our proof is represented visually in Figure \ref{proof diagram}.

\subsection{Geometry of the Markov numbers}
The question of counting for Markov triples,
 that is for $R>0$ estimating the size of the set of Markov triples
$(x,y,z)$ for which the sup norm is less than $R$, 
was first investigated by 
Gurwood  
in his thesis.
He  established an asymptotic formula using a correspondence
between Markov and Farey trees.
Somewhat later Zagier \cite{za} obtained an improved  error term. 
At the time of writing
the  best result   is due to a Rivin and the author \cite{mcr}:
$$M(R) = C(\log R)^2 + O(\log R \log \log R),$$
for a numerical constant which has a geometric interpretation.

\begin{figure}[hb]
\centering
\includegraphics[scale=.3]{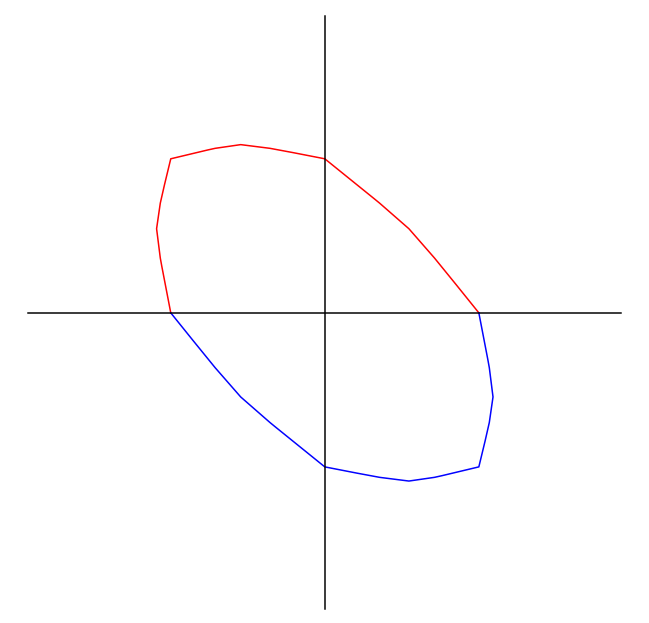} 
  \caption{A level set $S$ for the the length function $\ell$.}
  \label{pt graphic}
\end{figure}

The method introduced in \cite{mcr} differs from that of Gurwood and Zagier
in that it relies on  ideas of H. Cohn  to interpret the question 
in terms of geometry of hyperbolic surfaces.
In \cite{cohn} Cohn 
shows that the Markov triple $(1,1,1)$ is, 
up to multiplication by 3,
the image under the character map  $\chi$ of the 
modular torus.
Recall that the 
\textit{modular torus}
is the quotient of the upper half $\HH$
plane by  $\Gamma$,
the commutator subgroup of $\text{PSL}(2, \ZZ)$, acting by Mobius transformations.
Secondly, he shows that 
the permutations and the Vieta flips
used to construct Markov's binary tree
are induced by automorphisms of the
fundamental group of the torus
and concludes that  $\aut$
acts transitively on Markov triples.
It is well known that 
if $\rho : \F2 \rightarrow  \text{PSL}(2, \RR)$
is a discrete faithful representation and
$\gamma$  a closed curve  homotopic  
on the surface  $X = \HH/\rho(\F2)$,
then $\ell_X(\gamma)$, the length of the closed geodesic  homotopic to $\gamma$ 
on the surface  $X$,
can be computed from the trace using the relation
$$2\cosh \ell_X(\gamma)/2 = |\tr \rho(\gamma)|.$$
We note that the absolute value 
of the trace of an element 
of $\text{PSL}(2, \RR)$ is well defined.
Now identifying $\F2$ with the fundamental group of the modular torus
Cohn shows that  
if $x$ is a  Markov number then  there exists
a simple loop $\gamma$
such that 
$$x = \frac13 |\tr \rho(\gamma)|.$$
so there is an explicit relation between Markov numbers and
lengths of simple geodesics on the modular torus.
In  \cite{mcr} the simple closed curves are embedded 
as a subset of primitive elements in the integer lattice $\ZZ^2 \subset \RR^2$
and the length function $\gamma \mapsto \ell_X(\gamma)$ 
is shown, 
using a simple geometric argument,
 to extend to a proper  convex function $\ell: \RR^2 \mapsto \RR_+$.
 Thus the counting problem for Markov numbers 
 can be settled by counting integer points in the set
 $$\ell^{-1}([0,R])  = R\ell^{-1}([0,1]) \subset \RR^2$$
 It is easy to visualize this set
 using a very short computer program
 see for example Figure \ref{pt graphic}.
 
 In summary, 
 given a fraction $p/q$ 
 one obtains a primitive 
 lattice point $(q,p) \in \ZZ^2$
 and so a closed simple geodesic
 $\gamma_{p/q} \subset X$ 
 such that
 $$ M_{p/q} = 
 \frac23 \cosh\left(\frac{\ell_X(\gamma_{p/q}))}{2}\right)
 = 
 \frac23 \cosh(\| (q,p) \|_s),$$
 where $\|.\|_s$  is the norm 
 induced by the function
 $\frac12 \ell: \RR^2 \mapsto \RR_+$.
 
\subsection{Restatement of Conjectures}
Thus the fixed numerator conjecture
is equivalent to : 
Let $p, q$ and $i$ be positive integers
such that $i > 0, p < q$
and $p$ is coprime with 
both  $q, q+i$  then
$$\|(q,p) \|_s < \|(q + i,p) \|_s .$$
It is natural now to ask oneself
whether we need to restrict to integers.
In fact, using quite simple geometric arguments one can prove a result,
in a non discrete setting,
which yields Aigner's conjectures 
as a corollary:

\begin{thm}\label{main}
Let $p, q$
be real non negative numbers
and $i > 0$
then

\begin{equation}
\|(q,p) \|_s < \|(q + i,p) \|_s 
\end{equation}

and 

\begin{equation}
\|(q,p) \|_s < \|(q ,p +i ) \|_s 
\end{equation}

If in addition $p < q$ then

\begin{equation}
\|(q ,p  ) \|_s < \|(q + i ,p -i ) \|_s 
\end{equation}

\end{thm}

\subsection{Organisation, Remarks}

We give a brief account of 
the geometry and topology 
of the once punctured torus 
in Sections 2 and 3.
This is followed in Section 4 by the proof
of the convexity of the
length function $\ell$
which appeared in \cite{mcr}.
We conclude with the proof of 
Theorem \ref{main} in Section 5.

One could of course 
state many more theorems using 
these ideas but our purpose 
here is simply to show why
Aigner's conjectures are manifestations
of the convexity of the length function.

Finally we note that it is unlikely 
that purely geometric methods would
yield a proof of Frobenius' conjecture
as a natural variation of it is false \cite{mcp}.

\subsection{Thanks}

I thank Vlad Segesciu and Louis Funar for  many useful conversations over the years concerning this subject.
I am grateful too to Xu Binbin  for comments on the preliminary versions of the text.

\section{Hyperbolic structures and representations}

The rank two free group $\F2$ 
arises as the fundamental group of 
exactly two non-homeomorphic orientable surfaces:
the three-holed sphere $\Sigma_{0,3}$, 
the one-holed torus $\tt$, 
Furthermore $\tt$ enjoys the remarkable property
that every automorphism of its fundamental group
 is \textit{geometric} ie induced by a homeomorphism. 
Equivalently, every homotopy-equivalence $\tt \rightarrow \tt$
is homotopic to a homeomorphism. 
Much time and energy  has been 
devoted  to the geometric and topological 
properties of the  one holed torus
see for example \cite{buser}. 

\subsection{Fricke theory}

It  has been known since the time of Fricke that,
as a consequence of Riemann's Uniformization Theorem, 
the moduli space of hyperbolic structures
on the three holed sphere 
and the one holed torus 
can be identified with semi algebraic subsets of $\RR^3$. 
More precisely, 
one obtains a  hyperbolic structure  $X$ on $\Sigma$
from a discrete faithful representation $\rho$ of $\pi_1(\Sigma)$ 
 into $\isom$ such that $\HH/\rho(\pi_1(\Sigma))$  is homeomorphic to $\Sigma$.
One can lift any such  $\rho$
 from a representation into  $\isom^+$,
  which is isomorphic to  $\text{P}\sl2$, 
  to a representation
 $$\tilde{\rho} : \pi_1(\Sigma) \rightarrow \text{SL}(2,\RR).$$
Doing this allows one to consider the trace  $\tr \tilde{\rho}(\gamma) \in \RR$ 
for any element $\gamma$ of the free group 
so that many problems can be reduced to questions of  linear algebra.
In particular one defines a character map 
$$\chi :  X \mapsto (\tr \tilde{\rho} (\alpha) , \tr \tilde{\rho} (\beta), 
\tr \tilde{\rho} (\alpha\beta)),$$
which provides an embedding of the moduli space of hyperbolic structures on $\Sigma$.
Goldman \cite{gold} 
studied the dynamical system defined by the 
action of the group of outer automorphisms 
of $\pi_1(\Sigma)\simeq \F2$ on the moduli space.
If $\phi$ is an automorphism then 
it acts on the representations 
$$ \tilde{\rho} \mapsto \tilde{\rho}  \circ \phi^{-1}$$
and evidently this induces an action 
on the points in the image of the embedding above.
It turns out that,
 by applying  Cayley-Hamilton theorem to $\sl2$,
it is easy to show that this action
is the restriction of polynomial diffeomorphisms  of $\RR^3$.

\subsection{Markov triples}
An important aspect of the theory of the 
action of the group of outer automorphisms 
is that it admits an invariant function.
In the case of the pair of orientable surfaces this is 
$$ \kappa :  (x,y,z) \mapsto x^2 + y^2 + z^2 - x y z.$$
The level set $\kappa^{-1}(0)$ 
has five connected components - 
the singleton $(0,0,0)$
 and four codimension one manifolds homeomorphic to $\RR^2$.
 The non trivial components can be identified with the 
 Teichmueller space of a once punctured torus.
 It is well known that there are  integer points in $\kappa^{-1}(0) \cap \RR_+^3$,
 for example $(3,3,3)$,
 and that they form a single orbit for the action of the outer automorphisms.
 From such a  point one obtains a \text{Markov triple}, 
 that is a solution in positive integers  of the Markov cubic 
 $$x^2 + y^2 + z^2 - 3x y z = 0 $$
 simply by dividing through by $3$.
 This map is \textit{natural} in that it is a conjugation
 of automorphisms of the level set $\kappa^{-1}(0)$
 with the automorphisms of the Markov cubic.
So, 
since we will not be interested in arithemetic properties of these solutions,
we will also refer in what follows to elements of
 $\kappa^{-1}(0) \cap \NN^3 $ as Markov numbers.
 An integer  which appears in a Markov triple is called
 a \textit{Markov number}.
 
 \subsubsection{Reduction theory}
 
 The fact that the Markov triples form 
 a single orbit is a corollary of a classical result of Markov 
 in reduction theory which we now explain briefly.
 We start by giving the Markov triples an ordering,
 in the obvious way,
 using the sup norm on $\RR^3$ and proceed to show that
 any solution can be obtained from $(1,1,1)$
 by repeatedly applying generators of the  group of automorphisms 
 of the cubic.
This  automorphism group can be shown to be generated by:
 \begin{enumerate}
 \item sign change automorphisms $(x,y,z) \mapsto (x,-y,-z)$.
 \item coordinate permutations eg $(x,y,z) \mapsto (y,x,z)$.
 \item a Vieta flip $(x,y,z) \mapsto (x,y,3xy - z)$.
 \end{enumerate}
 Since Markov triples are solutions in positive integers
 we will use just
 the \textit{Markov morphisms} 
that is  the group generated 
 by permutations and the  Vieta flip.
 Given a Markov triple 
 $(x,y,z) \neq (1,1,1)$  we may apply a permutation 
 so that  $x \leq y < z$
 which one can try to ``reduce" using the Vieta flip, 
 that is replacing the it  by the triple  $(x,y,3xy - z)$.
 We obtain a smaller solution provided $3xy < 2z$
 and this  inequality  holds 
 for every Markov triple \textit{except}
 $(1,1,1)$.
 The reduction process gives rise
 to the structure of a rooted binary tree,
 the  \textit{Markov tree},
 on the Markov triples with $(1,1,1)$ as the \textit{root triple}.

\section{Automorphisms of $\F2$}
 
Let $\F2$ denote the free group on 2 generators
which we will denote $\alpha$ and $\beta$.
We think of  $\alpha$ and $\beta$ as simple loops 
meeting in a single point in a holed torus.

\subsection{Primitive elements}
An element $\gamma \in \F2$ is called \textit{primitive}, 
if  there exists an automorphism $\phi$, such that 
$\gamma = \phi(\alpha)$. 
If $\phi(\delta) = \beta$, then $\gamma$ and $\delta$ are called
\textit{associated primitives}.

Following \cite{mcr}
let $\phi : \F2 \rightarrow \ZZ^2$
be the canonical abelianizing homomorphism.
The kernel of $\phi$ is a \textit{characteristic subgroup}, 
that is it is $\aut$-invariant, 
so there is a homorphism 
from $\aut$ to the automorphisms of $\ZZ^2$
namely $\gl2$.
In fact this homomorphism is surjective and  faithful. Moreover,

\begin{enumerate}
\item
If $\gamma$ is primitive then $\phi(\gamma)$ 
is a primitive element of $\ZZ^2$.
\item
$\aut$ acts transitively on primitive elements of $\F2$
\item $\gl2$ acts transitively on primitive elements of $\ZZ^2$.

\end{enumerate}

\subsection{Topological realizations of automorphisms of $\F2$}

It is useful,
and this is essentially Cohn's point of view
\cite{cohn},
to identify $\F2$ with the fundamental group of  the holed torus
and  think of $\mathrm{Out}(\F2)$, the outer automorphism group of $\F2$,
as being  the mapping class group of the holed torus.
The outer automorphism group is generated 
by the so-called \textit{Nielsen transformations}, which either permute the basis 
$\alpha, \beta$, or transform it into $\alpha \beta, \beta$, or
$\alpha \beta^{-1}, \beta$. 
Both these  latter transformations can be realized  topologically as  Dehn twists on the holed torus.
We identify $\ZZ^2$ with the first homology of the holed torus $\Sigma$
and  will call the cover $\hc \rightarrow \Sigma$ corresponding 
to $\ker \phi$ the \textit{homology cover or maximal abelian cover of $\Sigma$}.


\section{Convexity of the length function}

Let $\gamma \in \F2$ be some closed curve 
and $\ell_X(\gamma)$ 
the minimum over the lengths  wrt the hyperbolic 
structure $X$ of all closed 
curves  in the homotopy 
class of $\gamma$.
In fact the  minimum is attained and is 
lthe length of the unique (oriented)  closed geodesic
in the homotopy  class.

\begin{thm}
The  length function is convex on $H_1(X,\QQ)$
and so extends continuously (as a convex function) to $H_1(X,\RR)$.
\end{thm}

Since the proof is short and instructive we shall reproduce it now.

Let $X$ be a holed torus equipped with a hyperbolic
structure. 
If $h$ is a primitive homology class 
(that is, not a multiple of another class), then the shortest multicurve representing a non-trivial homology
class $h$ is a  simple closed geodesic   
and a multiple of such a  geodesic otherwise. 
Further, the shortest multicurve representing $h$ is unique.

Begin by defining a valuation $\ell$ on the first homology  with integral coefficients, 
where $\ell(h)$ is just the length of the shortest multicurve representing h. 
The valuation of the trivial homology class is defined to be 0. 
By the previous paragraph $\ell$ statisfies
$$ \ell(nh) = n\ell(h), \forall n \in \NN.$$
Moreover if $h,g$ are not commensurable 
then $\ell$  satisfies the strict triangle inequality 
\begin{equation*} \label{triangle inequality}
\ell(h + g) < \ell(h) + \ell(g).
\end{equation*}
This is because the union of the shortest
multicurves representing  $h$ and $g$ is not embedded,
that is there is a transverse intersection somewhere.
Hence, the shortest multicurve representing 
$h + g$
 is strictly shorter than the union
 of the minimisers for $h,g$.
 
\section{Proof of Theorem \ref{main}.}

\begin{figure}[ht]
\centering
\includegraphics[scale=.3]{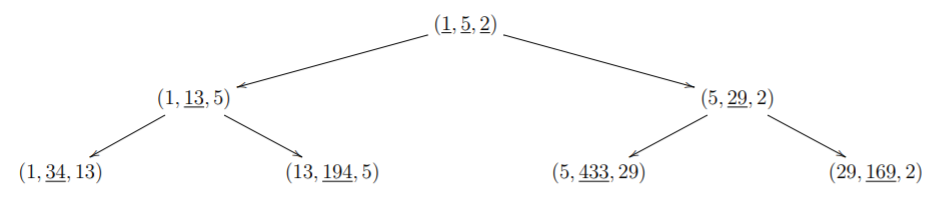} 
  \caption{Triples in a subtree of Markov's tree.}
  \label{markov tree}
\end{figure}

\begin{figure}[ht]
\centering
\includegraphics[scale=.3]{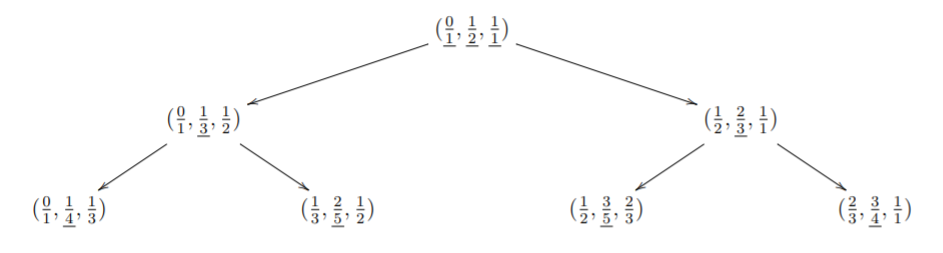} 
  \caption{Corresponding rationals in the Farey tree.}
  \label{farey tree}
\end{figure}

The correspondence 
between rationals and Markov
numbers is determined
by the values at
$0/1, 1/2, 1/1$
and,
following \cite{ra} 
(see Figures \ref{markov tree}
and \ref{farey tree}),
these are by convention $(1,5,2)$.
Now, passing from these fractions
to primitive  lattice points
we have the values of the norm at
$(1,0), (2,1), (1,1)$
We can easily find the value at
$(0,1)$ using Vieta flipping
and we see that it is $1$
ie the same value as for $(1,0)$.
In fact one can go a step further
and show that the three vectors
$(1,0), (0,1), (-1,1)$
correspond 
to the fundamental triple of Markov numbers
$(1,1,1)$.
The unit ball of the associated  stable norm looks just like
 the convex set depicted in Figure \ref{pt graphic}.

Now let $p, q$ 
be non negative
real numbers as in the statement of 
Theorem \ref{main}.
For any $i$ 
the point $(q + i, p)$ 
lies on the line $y = p$.
Since 
 $\ell: \RR^2 \mapsto \RR_+$
is a proper  convex function
this line meets
the set
$$S = \{ (x,y)\in \RR^2,\, 
\ell(x,y) = \ell(q,p) \} $$
in exactly two points
one of which is $(q,p)$
and the other  has
a negative $x$ coordinate.
It follows that for $i>0$
the point $(q + i, p)$ 
is outside 
of the norm ball 
$$B = \{ (x,y)\in \RR^2,\, 
\ell(x,y) \leq \ell(q,p) \} $$
and the first part of the theorem follows.

A similar argument using 
a vertical line instead of a
horizontal one yields a proof 
of the second part of the Theorem \ref{main}.

Finally for the third part 
one must consider the line
$(q,p) + i(-1,1)$.
It is easy to check that 
it meets $S$ in $(q,p)$ 
and $(p,q)$ and so for
$i >0$ we have the required inequality.

%

\end{document}